\pdfoutput=1
\documentclass[pdflatex,preprint,number]{elsarticle}

\usepackage{amssymb}
\usepackage{amsfonts}
\usepackage{hyperref}
\usepackage{epsfig,amsmath}
\usepackage{tikz,ifthen,calc}
\usepackage{graphicx}
\usepackage{caption}
\usepackage{subcaption}
\usepackage{wrapfig}

\usepackage[margin=1.00 in]{geometry}

\newtheorem{thm}{Theorem}
\newtheorem{lem}[thm]{Lemma}
\newtheorem{coro}[thm]{Corollary}
\newdefinition{rmk}{Remark}
\newdefinition{definition}{Definition}
\newdefinition{example}{Example}
\newproof{pf}{Proof}
\newproof{pot}{Proof of Theorem \ref{thm2}}

\title{A rigidity property of ribbon $L$-shaped $n$-ominoes and generalizations\tnoteref{t1}}
\tnotetext[t1]{V. Nitica was partially supported by Simons Foundation Grant 208729.}

\author[vn]{Viorel Nitica}
\ead{vnitica@wcupa.edu}

\address[vn]{Department of Mathematics, West Chester University, West Chester, PA 19380\\ and Institute of Mathematics, Romanian Academy, P.O. Box 1--764, RO-70700 Bucharest, Romania}

\begin{document}

\begin{abstract} Let $n\ge 4$ even and let $\mathcal{T}_n$ be the set of ribbon $L$-shaped $n$-ominoes. We study tiling problems for regions in a square lattice by $\mathcal{T}_n$. Our main result shows a remarkable rigidity property: a tiling of the first quadrant by $\mathcal{T}_n$ is possible if and only if it reduces to a tiling by $2\times n$ and $n\times 2$ rectangles. An application is the classification of all rectangles that can be tiled by $\mathcal{T}_n$: a rectangle can be tiled by $\mathcal{T}_n$ if and only if both of its sides are even and at least one side is divisible by $n$. Another application is the existence of the local move property for an infinite family of sets of tiles: $\mathcal{T}_n$ has the local move property for the class of rectangular regions with respect to the local moves that interchange a tiling of an $n\times n$ square by $n/2$ vertical  rectangles, with a tiling by $n/2$ horizontal rectangles, each vertical/horizontal rectangle being covered by two ribbon $L$-shaped $n$-ominoes. We show that these results are not valid for any $n$ odd. The rectangular pattern of a tiling persists if we add an extra $2\times 2$ tile to $\mathcal{T}_n$. A rectangle can be tiled by the larger set of tiles if and only if it has both sides even. In contrast, the addition of an extra even $\times$ odd or odd $\times$ odd rectangle to one of the above sets of tiles allows for a tiling of the first quadrant that does not respect the rectangular pattern.
\end{abstract}




\maketitle

\section{Introduction}\label{intro}

In this article we study tiling problems for regions in a square lattice by certain symmetries of an $L$-shaped polyomino.
Polyominoes were introduced by Solomon W. Golomb \cite{Golomb1} and the standard reference about this subject is the book \emph{Polyominoes} \cite{Golumb3}.
They are a never ending source of combinatorial problems.

The $L$-shaped polyomino we study is placed in a square lattice and is made out of $n, n\ge 3,$ unit squares, or \emph{cells}. See Figure~\ref{fig:LTetromino1}. In an $a\times b$ rectangle, $a$ is the height and $b$ is the base. We consider translations (only!) of the tiles shown in Figure~\ref{fig:LTetrominoes1}. The first four are ribbon $L$-shaped $n$-ominoes and the last one is a $2\times 2$ square. A {\em ribbon polyomino}~\cite{Pak} is a simply connected polyomino with no two unit squares lying along a line parallel to $y=x$. We denote the set of tiles $\{T_1, T_2, T_3, T_4, T_5\}$ by $\mathcal{T}_n^+$ and the set of tiles $\{T_1, T_2, T_3, T_4\}$ by $\mathcal{T}_n.$

\begin{figure}[h]
~~~~~~
\begin{subfigure}{.15\textwidth}
\begin{tikzpicture}[scale=.4]
\draw [line width = 1](2.5,1)--(1,1)--(1,2)--(0,2)--(0,0)--(2.5,0);
\draw[dotted, line width=1] (2.6,.5)--(3.9,.5);
\draw [line width = 1] (4,0)--(6,0)--(6,1)--(4,1);
\draw [line width = 1] (1,0)--(1,1);
\draw [line width = 1] (0,1)--(1,1);
\draw [line width = 1] (2,0)--(2,1);
\draw [line width = 1] (5,0)--(5,1);
\end{tikzpicture}
\caption{An $L$ $n$-omino with $n$ cells.}
\label{fig:LTetromino1}
\end{subfigure}
~~~~~~~~~
\begin{subfigure}{.5\textwidth}
\begin{tikzpicture}[scale=.4]
\draw [line width = 1](2.5,1)--(1,1)--(1,2)--(0,2)--(0,0)--(2.5,0);
\draw[dotted, line width=1] (2.6,.5)--(3.9,.5);
\draw [line width = 1] (4,0)--(6,0)--(6,1)--(4,1);
\draw [line width = 1] (1,0)--(1,1);
\draw [line width = 1] (0,1)--(1,1);
\draw [line width = 1] (2,0)--(2,1);
\draw [line width = 1] (5,0)--(5,1);
\node at (3, -1) {$T_3$};
\draw [line width = 1] (10,1)--(8,1)--(8,2)--(10,2);
\draw [dotted, line width=1] (10.1,1.5)--(11.4,1.5);
\draw [line width = 1] (9,1)--(9,2);
\draw [line width = 1] (11.5,1)--(13,1)--(13,0)--(14,0)--(14,2)--(11.5,2);
\draw [line width = 1] (12,1)--(12,2);
\draw [line width = 1] (13,1)--(13,2);
\draw [line width = 1] (13,1)--(14,1);
\node at (11, -1) {$T_4$};

\draw [line width = 1] (16,0)--(18,0)--(18,2)--(16,2)--(16,0);
\draw [line width = 1] (17,0)--(17,2);
\draw [line width = 1] (16,1)--(18,1);
\node at (17, -1) {$T_5$};

\draw [line width = 1] (-9,2)--(-9,4)--(-8,4)--(-8,2);
\draw [dotted, line width=1] (-8.5,1.9)--(-8.5,.6);
\draw [line width = 1] (-9,.5)--(-9,-2)--(-7,-2)--(-7,-1)--(-8,-1)--(-8,.5);
\draw [line width = 1] (-9,3)--(-8,3);
\draw [line width = 1] (-9,0)--(-8,0);
\draw [line width = 1] (-9,-1)--(-8,-1);
\draw [line width = 1] (-8,-1)--(-8,-2);
\node at (-8, -3) {$T_1$};

\draw [line width = 1] (-3,1.5)--(-3,3)--(-4,3)--(-4,4)--(-2,4)--(-2,1.5);
\draw [line width = 1] (-3,2)--(-2,2);
\draw [line width = 1] (-3,3)--(-2,3);
\draw [line width = 1] (-3,3)--(-3,4);
\draw [dotted, line width=1] (-2.5,1.4)--(-2.5,.1);
\draw [line width = 1] (-3,0)--(-3,-2)--(-2,-2)--(-2,0);
\draw [line width = 1] (-3,-1)--(-2,-1);
\node at (-2.5, -3) {$T_2$};
\end{tikzpicture}
\caption{The set of tiles $\mathcal{T}^+$.}
\label{fig:LTetrominoes1}
\end{subfigure}
\caption{}
\end{figure}

The inspiration for this paper is the recent publication~\cite{CLNS}, showing related results for the set of tiles consists of ribbon $L$-tetrominoes. That is, \cite{CLNS} investigates tiling problems for the set of tiles $\mathcal{T}_n$ in the particular case $n=4$. The extension of the results in~\cite{CLNS} to the case of  general ribbon $L$-shaped $n$-ominoes is a natural question. A coloring invariant makes the proofs in \cite{CLNS} more transparent then here. The invariant does not generalize to this new situation and more involved geometric arguments are developed in this paper.

In order to avoid repetition, we assume for the rest of the paper that, unless otherwise specified, $n$ is even and $n\ge 4$. We denote the first quadrant by $Q_1$.

\begin{definition} A tiling of $Q_1$ by $\mathcal{T}_n^+$ {\em follows the rectangular pattern}  if it reduces to a tiling by $2\times 2, 2\times n, n\times 2$ rectangles, with the last two types being covered by two ribbon $L$-shaped $n$-ominoes. More general, let $P$ be a polygonal region in $Q_1$. Then $P$ {\em follows the rectangular pattern} in a tiling of $Q_1$ by $\mathcal{T}_n^+$ if $P$ is completely covered by non overlapping $2\times 2, 2\times n$ and $n\times 2$ rectangles having the coordinates of all vertices even and with each $2\times n$ or $n\times 2$ rectangle covered by two ribbon $L$-shaped $n$-ominoes. Similar notions can be introduced for $\mathcal{T}_n$.
\end{definition}

Our main result is the following.

\begin{thm}\label{main1} Any tiling of $Q_1$ by $\mathcal{T}_n^+$ or $\mathcal{T}_n$ follows the rectangular pattern. Moreover, the addition of an even $\times$ odd or odd $\times$ odd rectangle to the set of tiles $\mathcal{T}_n^+$ or $\mathcal{T}_n$ allows for a tiling of $Q_1$ by $\mathcal{T}_n^+$ (respectively $\mathcal{T}_n$) that does not follow the rectangular pattern.
\end{thm}

Some immediate consequences of the main result are listed below. The proofs of Corollaries \ref{rectangle2}, \ref{inflated1}, \ref{thm-half-infinite-strip} are similar to those of \cite[Corollary 2, 3, 4, 5]{CLNS}.

\begin{coro}\label{rectangle2} A rectangle can be tiled by $\mathcal{T}_n^+$ if and only if the tiling follows the rectangular pattern. Consequently, a rectangle can be tiled by $\mathcal{T}_n^+$ if and only if both of its sides are even.
\end{coro}

A $k$-copy of a polyomino is a replica of it in which all $1\times 1$ squares are replaced by $k\times k$ squares.

\begin{coro}\label{inflated1} Let $k$ odd. Then a $k$-copy of the ribbon $L$ $n$-omino cannot be tiled by $\mathcal{T}_n^+$.
\end{coro}

\begin{figure}[h]
\begin{subfigure}{.28\textwidth}
\begin{tikzpicture}[scale=.3]
\draw [line width = 1](0,0)--(9,0)--(9,1)--(1,1)--(1,2)--(0,2)--(0,0);
\draw [line width = 1] (1,1)--(3,1)--(3,2)--(2,2)--(2,10)--(1,10)--(1,1);
\draw [line width = 1] (3,1)--(5,1)--(5,2)--(4,2)--(4,10)--(3,10)--(3,1);
\draw [line width = 1] (5,1)--(5,10)--(6,10)--(6,2)--(7,2);
\draw [line width = 1] (3,2)--(3,11)--(1,11)--(1,10);
\draw [line width = 1] (5,2)--(5,11)--(3,11)--(3,10);
\draw [line width = 1] (7,2)--(7,11)--(5,11)--(5,10);
\draw [line width = 1] (9,2)--(9,11)--(7,11)--(7,10);
\draw [line width = 1] (1,11)--(1,12)--(10,12)--(10,10)--(9,10);
\draw [line width = 1] (7,1)--(9,1)--(9,2)--(8,2)--(8,10)--(7,10)--(7,1);
\draw [line width = 1] (-2,0)--(12,0);
\draw [line width = 1] (-2,12)--(12,12);
\draw [line width = 1] (0,2)--(-1,2)--(-1,12);
\draw [line width = 1] (-1,2)--(-2,2);
\draw [line width = 1] (10,10)--(11,10)--(11,0);
\draw [line width = 1] (11,10)--(12,10);
\end{tikzpicture}
\caption{An infinite strip.}
\label{fig:InfiniteTiling-halfplane-bis}
\end{subfigure}
~
\begin{subfigure}{.35\textwidth}
		\begin{tikzpicture}[scale=.3]
			\draw[dotted] (-16.5, 0) grid (0, 14.5);
			\draw[line width = 1, ->] (0, 0) -- (-16.5, 0);
			\draw[line width = 1, ->] (0, 0) -- (0, 14.5);
			\draw[line width = 1] (-5,0)--(-5,2)--(-4,2)--(-4,1)--(0,1);
			\draw[line width = 1] (-2,1)--(-2,7)--(0,7);
            \draw[line width = 1] (-2,7)--(-4,7)--(-4,2);
            \draw[line width = 1] (-2,7)--(-2,13)--(0,13);
            \draw[line width = 1] (-2,13)--(-4,13)--(-4,7);
            \draw[line width = 1] (-2,13)--(-2,14.5);
			\draw[line width = 1] (-4,13)--(-4,14.5);
            \draw[line width = 1] (-5,2)--(-11,2)--(-11,0);
            \draw[line width = 1] (-11,2)--(-16,2);
            \draw[line width = 1] (-6,2)--(-6,8)--(-4,8);
            \draw[line width = 1] (-6,8)--(-6,14)--(-4,14);
            \draw[line width = 1] (-8,2)--(-8,8)--(-6,8);
            \draw[line width = 1] (-8,8)--(-8,14)--(-6,14);
            \draw[line width = 1] (-10,2)--(-10,8)--(-8,8);
            \draw[line width = 1] (-10,8)--(-10,14)--(-8,14);
            \draw[line width = 1] (-12,2)--(-12,8)--(-10,8);
            \draw[line width = 1] (-12,8)--(-12,14)--(-10,14);
            \draw[line width = 1] (-14,2)--(-14,8)--(-12,8);
            \draw[line width = 1] (-14,8)--(-14,14)--(-12,14);
            \draw[line width = 1] (-16,2)--(-16,8)--(-14,8);
            \draw[line width = 1] (-16,8)--(-16,14)--(-14,14);
            \draw [line width = 1] (-16,2)--(-16.5,2);
            \draw [line width = 1] (-16,8)--(-16.5,8);
            \draw [line width = 1] (-16,14)--(-16.5,14);
            \draw [line width = 1] (-16,14)--(-16,14.5);
            \draw [line width = 1] (-14,14)--(-14,14.5);
            \draw [line width = 1] (-12,14)--(-12,14.5);
            \draw [line width = 1] (-10,14)--(-10,14.5);
            \draw [line width = 1] (-8,14)--(-8,14.5);
            \draw [line width = 1] (-6,14)--(-6,14.5);
		\end{tikzpicture}
\caption{The second quadrant.}
\label{fig:second_quadrant-bis}
\end{subfigure}
~
\begin{subfigure}{.28\textwidth}
\begin{tikzpicture}[scale=.3]
\draw [line width = 1.5] (0,0)--(16,0)--(16,15)--(0,15)--(0,0);
\draw [line width = 1.5] (0,10)--(5,10)--(5,0);
\draw [line width = 1.5] (5,10)--(6,10)--(6,15);
\draw [line width = 1.5] (6,10)--(6,9)--(9,9)--(9,8)--(5,8)--(5,10);
\draw [line width = 1.5] (9,9)--(11,9)--(11,5)--(10,5)--(10,8)--(9,8);
\draw [line width = 1.5] (10,0)--(10,5)--(16,5);
\draw [line width = 1.5] (11,9)--(11,15);
\node at (13, 2.5) {\small{$(n,n+1)$}};
\node at (3, 12.5) {\small{$(n,n+1)$}};
\node at (2.5, 5) {\small{$(2n,n)$}};
\node at (13.5, 10) {\small{$(2n,n)$}};
\node at (7.5, 4) {\small{$(2n-2,n)$}};
\node at (8.4, 11.5) {\small{$(n+1,n)$}};
\end{tikzpicture}
\caption{A $(3n,3n+1)$ rectangle.}
\label{fig:odd-rectangle}
\end{subfigure}
\caption{}
\end{figure}

\begin{coro}\label{thm-half-infinite-strip}
 A half-infinite strip of odd width cannot be tiled by $\mathcal{T}_n^+$.
\end{coro}

It was proved by de~Brujin~\cite{Brujin} that an $a\times b$ rectangle can be tiled by $k\times 1$ and $1\times k$ bars if and only if $k$ divides one of the sides of the rectangle. In conjunction with Theorem~\ref{main1} this gives:

\begin{thm}\label{rectangle1} A rectangle can be tiled by $\mathcal{T}_n$ if and only if both sides are even and at least one of the sides is divisible by $n$. Any tiling of a rectangle by $\mathcal{T}_n$ has to follow the rectangular pattern.
\end{thm}

\begin{rmk} Not much is known about tiling integer size rectangles with an odd size if we allow in the set of tiles all 8 orientations of an $L$-shaped $n$-omino, $n$ even. It is known that if $n\equiv 0$ mod 4, then the area of a rectangle that can be tiled is a multiple of $2n$. This condition is sufficient if $n=4$ and if the rectangle is not a bar of height 1, see for example~\cite{Nit-bis}, but not in general.
\end{rmk}

\begin{rmk} Figure~\ref{fig:InfiniteTiling-halfplane-bis} shows a tiling by $\mathcal{T}_n$ of an infinite strip of width $n+2$ that does not follow the rectangular pattern. Figure~\ref{fig:second_quadrant-bis} shows a tiling of the second quadrant by $\mathcal{T}_n$ that does not follow the rectangular pattern. In both figures, the small rectangles are tiled by two ribbon $L$-shaped $n$-ominoes. The last example shows that our results do not remain valid for the reflections of the tiled region about the horizontal/vertical axis.
\end{rmk}

\begin{rmk} Assume $n$ odd. Figure~\ref{fig:odd-rectangle} shows how to tile a $(3n)\times (3n+1)$ rectangle by $\mathcal{T}_n$. Each interior rectangle in Figure~\ref{fig:odd-rectangle} can be covered by $2\times n$ or $n\times 2$ rectangles, which in turn are covered by two ribbon $L$-shaped $n$-ominoes. Thus the assumption that $n$ is even is necessary in Theorems~\ref{main1}, \ref{rectangle1}, and Corollaries~\ref{rectangle2},  \ref{inflated1}, \ref{thm-half-infinite-strip}. Due to the example in Figure~\ref{fig:odd-rectangle}, tiling of certain half strips of odd width by $\mathcal{T}_n$, is possible for any $n$ odd.
\end{rmk}

\begin{definition} Given a set of tiles $\mathcal{T}$ and a finite set of local replacement moves $\mathcal{L}$ for the tiles in $\mathcal{T}$, we say that a region $\Gamma$ has \emph{local connectivity} with respect to $\mathcal{T}$ and $\mathcal{L}$ if it possible to convert any tiling of $\Gamma$ into any other by means of these moves. If $\mathcal{R}$ is a class of regions, then we say that there is a \emph{local move property} for $\mathcal{T}$ and $\mathcal{R}$ if there exists a finite set of moves $\mathcal{L}$ such that every $\Gamma$ in $\mathcal{R}$ has local connectivity with respect to $\mathcal{T}$ and $\mathcal{L}$.
\end{definition}

It is shown in~\cite{KK} that if the set of tiles consists of $1\times k$ and $k\times 1$ bars, then the class of rectangular regions has the local move property with respect to the moves that interchange a tiling of an $k\times k$ square by $k$ vertical bars with a tiling of the same $k\times k$ square by $k$ horizontal bars. In conjunction with Theorem \ref{main1}, this shows:

\begin{thm}  The set of tiles $\mathcal{T}_n$ has the local move property for the class of rectangular regions. The local moves interchange an $n\times n$ square tiled by $n/2$ vertical $n\times 2$ rectangles with the same $n\times n$ square tiled by $n/2$ horizontal $2\times n$ rectangles. Each vertical/horizontal rectangle is covered by two ribbon $L$-shaped $n$-ominoes.
\end{thm}

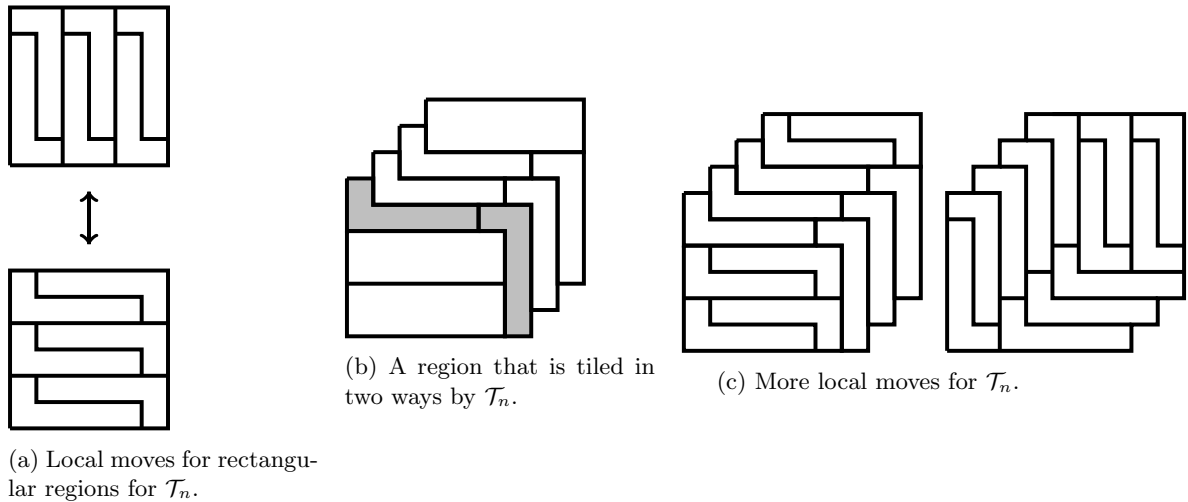
\begin{figure}[h]
\begin{subfigure}{.25\textwidth}
\begin{tikzpicture}[scale=.35]
\draw [line width = 1.5] (0,4)--(6,4)--(6,10)--(0,10)--(0,4);
\draw [line width = 1.5] (2,4)--(2,10);
\draw [line width = 1.5] (4,4)--(4,10);
\draw [line width = 1.5] (2,5)--(1,5)--(1,9)--(0,9);
\draw [line width = 1.5] (4,5)--(3,5)--(3,9)--(2,9);
\draw [line width = 1.5] (6,5)--(5,5)--(5,9)--(4,9);

\draw [line width = 1.5,->] (3,3)--(3,1);
\draw [line width = 1.5,->] (3,1)--(3,3);

\draw [line width = 1.5] (0,-6)--(6,-6)--(6,0)--(0,0)--(0,-6);
\draw [line width = 1.5] (1,0)--(1,-1)--(5,-1)--(5,-2);
\draw [line width = 1.5] (1,-2)--(1,-3)--(5,-3)--(5,-4);
\draw [line width = 1.5] (1,-4)--(1,-5)--(5,-5)--(5,-6);
\draw [line width = 1.5] (0,-2)--(6,-2);
\draw [line width = 1.5] (0,-4)--(6,-4);
\end{tikzpicture}
\caption{Local moves for rectangular regions for $\mathcal{T}_n$.}
\label{fig:local-moves-757}
\end{subfigure}
~
\begin{subfigure}{.25\textwidth}
\begin{tikzpicture}[scale=.35]
\draw [line width = 0] [fill = lightgray] (0,2)--(0,0)--(5,0)--(5,1)--(1,1)--(1,2)--(0,2);
\draw [line width = 0] [fill = lightgray] (5,0)--(5,1)--(7,1)--(7,-4)--(6,-4)--(6,0)--(5,0);
\draw [line width = 1.5] (0,2)--(0,0)--(5,0)--(5,1)--(1,1)--(1,2)--(0,2);
\draw [line width = 1.5] (5,0)--(5,1)--(7,1)--(7,-4)--(6,-4)--(6,0)--(5,0);
\draw [line width = 1.5] (1,3)--(1,1)--(6,1)--(6,2)--(2,2)--(2,3)--(1,3);
\draw [line width = 1.5] (6,1)--(6,2)--(8,2)--(8,-3)--(7,-3)--(7,1)--(6,1);
\draw [line width = 1.5] (2,4)--(2,2)--(7,2)--(7,3)--(3,3)--(3,4)--(2,4);
\draw [line width = 1.5] (7,2)--(7,3)--(9,3)--(9,-2)--(8,-2)--(8,2)--(7,2);
\draw [line width = 1.5] (3,5)--(3,3)--(9,3)--(9,5)--(3,5);
\draw [line width = 1.5] (0,0)--(6,0)--(6,-4)--(0,-4)--(0,0);
\draw [line width = 1.5] (0,-2)--(6,-2);
\end{tikzpicture}
\caption{A region that is tiled in two ways by $\mathcal{T}_n$.}
\label{fig:local-moves}
\end{subfigure}
~
\begin{subfigure}{.3\textwidth}
\begin{tikzpicture}[scale=.35]
\draw [line width = 1.5] (0,2)--(0,0)--(5,0)--(5,1)--(1,1)--(1,2)--(0,2);
\draw [line width = 1.5] (5,0)--(5,1)--(7,1)--(7,-4)--(6,-4)--(6,0)--(5,0);
\draw [line width = 1.5] (1,3)--(1,1)--(6,1)--(6,2)--(2,2)--(2,3)--(1,3);
\draw [line width = 1.5] (6,1)--(6,2)--(8,2)--(8,-3)--(7,-3)--(7,1)--(6,1);
\draw [line width = 1.5] (2,4)--(2,2)--(7,2)--(7,3)--(3,3)--(3,4)--(2,4);
\draw [line width = 1.5] (7,2)--(7,3)--(9,3)--(9,-2)--(8,-2)--(8,2)--(7,2);
\draw [line width = 1.5] (3,5)--(3,3)--(9,3)--(9,5)--(3,5);
\draw [line width = 1.5] (0,0)--(6,0)--(6,-4)--(0,-4)--(0,0);
\draw [line width = 1.5] (0,-2)--(6,-2);
\draw [line width = 1.5] (1,0)--(1,-1)--(5,-1)--(5,-2);
\draw [line width = 1.5] (1,-2)--(1,-3)--(5,-3)--(5,-4);
\draw [line width = 1.5] (4,5)--(4,4)--(8,4)--(8,3);

\draw [line width = 1.5] (10,-4)--(16,-4)--(17,-4)--(17,-3)--(18,-3)--(18,-2)--(19,-2)--(19,5)--(13,5)--(13,4)--(12,4)--(12,3)--(11,3)--(11,2)--(10,2)--(10,-4);
\draw [line width = 1.5] (10,1)--(11,1)--(11,-3)--(12,-3);
\draw [line width = 1.5] (11,2)--(12,2)--(12,-4);
\draw [line width = 1.5] (12,3)--(13,3)--(13,-3)--(17,-3);
\draw [line width = 1.5] (12,-2)--(13,-2);
\draw [line width = 1.5] (13,4)--(14,4)--(14,-2)--(18,-2);
\draw [line width = 1.5] (13,-1)--(14,-1);
\draw [line width = 1.5] (14,5)--(15,5)--(15,-1)--(19,-1);
\draw [line width = 1.5] (14,0)--(15,0);
\draw [line width = 1.5] (17,5)--(17,-1);
\draw [line width = 1.5] (15,4)--(16,4)--(16,0)--(17,0);
\draw [line width = 1.5] (17,4)--(18,4)--(18,0)--(19,0);
\end{tikzpicture}
\caption{More local moves for $\mathcal{T}_n$.}
\label{fig:local-moves-bis}
\end{subfigure}
\caption{Local moves for $\mathcal{T}_n$.}
\end{figure}

\begin{rmk} Figure~\ref{fig:local-moves-757} shows the local moves for $n=6$. For general regions, the set $\mathcal{T}_n$ does not have the local move property. For each $n$ there exists an infinite family of row-convex regions that admit only two tilings by $\mathcal{T}_n$. Figure~\ref{fig:local-moves} shows one representative, for $n=6$. In general, there are $n/2$ rectangles of size $n\times 2$ in the figure, one at the top and $n/2-1$ at the bottom, and  each one tiled by two tiles from $\mathcal{T}_n$.  The infinite family of regions is obtained by introducing more copies of the gray subregion in the middle. For $n=6$ the two tilings are shown in Figure~\ref{fig:local-moves-bis}.
\end{rmk}

\begin{rmk} One may wonder if $\mathcal{T}_n$ remains rigid if other tiles besides rectangles are added to the tiling set. Figure~\ref{final-ilustration-39} shows a tiling of a $6\times 10$ rectangle by $\mathcal{T}_4\cup \mathcal{T}_6$ that does not follow the rectangular pattern. The regions I, II are $4\times 4$ squares and can be tiled by $\mathcal{T}_4$. It is easy to see that this example can be generalized to the set of tiles $\mathcal{T}_m\cup \mathcal{T}_n$ for any $m\not =n$ with $m,n$ even.
\end{rmk}
		
\section{Tiling $Q_1$ by $\mathcal{T}_n^+$}\label{s:tilingbyTn}

In this section we prove Theorem~\ref{main1}. For simplicity we refer to the $2\times 2$ squares with the coordinates of all vertices even as \emph{even $2\times 2$ squares.}

\begin{definition}\label{d:def2} A $T_2$ tile that is part of a tiling of $Q_1$ by $\mathcal{T}_n^+$ is said to be \emph{in an irregular position} if the coordinates of its lowest left corner are even, and if all even $2\times 2$ squares below and to the left of its lowest left corner follow the rectangular pattern. The corresponding notions for a $T_4$ tile are defined via a reflection about the line $x=y$.
\end{definition}

\begin{rmk} The $T_2$ tile in Figure~\ref{fig:aux-lemma-098} is in irregular position. All gray squares follows the rectangular pattern.
\end{rmk}

\begin{definition}\label{d:def3} Assume that $Q_1$ is tiled by $\mathcal{T}_n^+$. A \emph{gap} is a rectangular region of height $2$ in the square lattice, that has the coordinates of the lower left corner even. If its length is even, the gap is called \emph{even}, otherwise it is called \emph{odd}. We assume that the even $2\times 2$ squares on the left side of the gap and below the upper level of the gap as well as the even $2\times 2$ squares directly below the gap follow the rectangular pattern. An even gap has a right vertical side of height 2 and an odd gap has a right vertical side of height 1. The right vertical side of a gap cannot be covered by tiles from the tiling of $Q_1$.
\end{definition}

\begin{rmk} Pictures of gaps are shown in Figure~\ref{fig-bis-34}. The the dark gray regions follow the rectangular pattern. We do not assume anything about the tiling of the remaining white regions in $Q_1$ or the gaps.
\end{rmk}

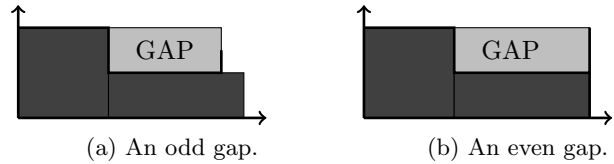
\begin{figure}[h]
\centering
\begin{subfigure}{0.25\textwidth}
\begin{tikzpicture}[scale=.3]

\draw[fill=lightgray] (4, 2) rectangle (9, 4);
\node at (6.5,3) {GAP};

\draw [fill=darkgray] (0,0) rectangle (4,4);
\draw [fill=darkgray] (4,0) rectangle (10,2) ;

\draw [line width = 1,->] (0,0)--(0,5);
\draw [line width = 1,->] (0,0)--(11,0);

\draw [line width = 1] (0,4)--(4,4)--(4,2)--(9,2)--(9,3);

\end{tikzpicture}
\caption{An odd gap.}
\label{fig-bis-34}
\end{subfigure}
~~
\begin{subfigure}{0.25\textwidth}
\begin{tikzpicture}[scale=.3]

\draw[fill=lightgray] (4, 2) rectangle (10, 4);
\node at (6.5,3) {GAP};

\draw [fill=darkgray] (0,0) rectangle (4,4);
\draw [fill=darkgray] (4,0) rectangle (10,2) ;

\draw [line width = 1,->] (0,0)--(0,5);
\draw [line width = 1,->] (0,0)--(11,0);

\draw [line width = 1] (0,4)--(4,4)--(4,2)--(10,2)--(10,4);

\draw [line width = 1] (10,2)--(10,0);
\end{tikzpicture}
\caption{An even gap.}
\label{fig-bis-34}
\end{subfigure}
\caption{Pictures of gaps.}
\end{figure}

For the following three lemmas we assume that a tiling of $Q_1$ by $\mathcal{T}_n^+$ is given.

\begin{lem}\label{lemma-123} Assume that the leftmost even $2\times 2$ square in an odd gap of length $L\ge 3$ that does not follow the rectangular pattern. Then there exists a $T_2$ tile in an irregular position that is above the bottom of the gap and to the left of the right side of the gap.
\end{lem}

\begin{figure}[h]
\center
\begin{subfigure}{0.2\textwidth}
		\centering
		\begin{tikzpicture}[scale=.35]
        \draw[fill=lightgray] (0,0) rectangle (5,2);
		\draw[line width = 1, ->] (0, -2) -- (7.5, -2);
		\node at (.5, .5) {$1$};
        \node at (1.5, 1.5) {$2$};
        \draw[fill=darkgray] (2,0)--(3,0)--(3,1)--(2,1);
		\draw[line width=1] (0, 0) -- (2, 0) -- (2, 1) -- (1, 1) -- (1, 5) -- (0, 5) -- (0, 0);
        \draw[line width=1] (1,1)--(1,6)--(2,6)--(2,2)--(3,2)--(3,1)--(2,1);
		\draw[dotted,step=1] (0, -2) grid (7.5, 6.5);
		\draw[line width = 1, ->] (0, -2) -- (0, 6.5);
		\end{tikzpicture}
		\caption{}
	\end{subfigure}
~
\begin{subfigure}{0.2\textwidth}
		\centering
		\begin{tikzpicture}[scale=.35]
        \draw[fill=lightgray] (0,0) rectangle (7,2);
		\draw[line width = 1, ->] (0, -2) -- (7.5, -2);
		\node at (.5, .5) {$1$};
        \node at (1.5, 1.5) {$2$};
        \draw[fill=darkgray] (2,0)--(3,0)--(3,1)--(2,1);
		\draw[line width=1] (0, 0) -- (2, 0) -- (2, 1) -- (1, 1) -- (1, 5) -- (0, 5) -- (0, 0);
        \draw[line width=1] (1,3)--(2,3)--(2,2)--(6,2)--(6,1)--(2,1);
		\draw[dotted,step=1] (0, -2) grid (7.5, 6.5);
		\draw[line width = 1, ->] (0, -2) -- (0, 6.5);
		\end{tikzpicture}
		\caption{}
	\end{subfigure}
~
\begin{subfigure}{0.2\textwidth}
		\centering
		\begin{tikzpicture}[scale=.35]
        \draw[fill=lightgray] (0,0) rectangle (7,2);
		\draw[line width = 1, ->] (0, -2) -- (7.5, -2);
		\node at (.5, .5) {$1$};
        \node at (1.5, 1.5) {$2$};
        \draw[fill=darkgray] (2,0)--(3,0)--(3,1)--(2,1);
		\draw[line width=1] (0, 0) -- (2, 0) -- (2, 1) -- (1, 1) -- (1, 5) -- (0, 5) -- (0, 0);
        \draw[line width=1] (1,2)--(6,2)--(6,0)--(5,0)--(5,1)--(1,1);
		\draw[dotted,step=1] (0, -2) grid (7.5, 6.5);
		\draw[line width = 1, ->] (0, -2) -- (0, 6.5);
		\end{tikzpicture}
		\caption{}
	\end{subfigure}
~
\begin{subfigure}{0.2\textwidth}
		\centering
		\begin{tikzpicture}[scale=.35]
        \draw[fill=lightgray] (0,0) rectangle (7,2);
		\draw[line width = 1, ->] (0, -2) -- (7.5, -2);
		\node at (.5, .5) {$1$};
        \node at (1.5, 1.5) {$2$};
        \draw[fill=darkgray] (2,0)--(3,0)--(3,1)--(2,1);
		\draw[line width=1] (0, 0) -- (2, 0) -- (2, 1) -- (1, 1) -- (1, 5) -- (0, 5) -- (0, 0);
        \draw[line width=1] (1,1)--(3,1)--(3,3)--(1,3)--(1,1);
		\draw[dotted,step=1] (0, -2) grid (7.5, 6.5);
		\draw[line width = 1, ->] (0, -2) -- (0, 6.5);
		\end{tikzpicture}
		\caption{}
	\end{subfigure}
\caption{A $T_1$ tile covers the lower leftmost cell in the gap-the base case.}
\label{fig-bis-35}
\end{figure}
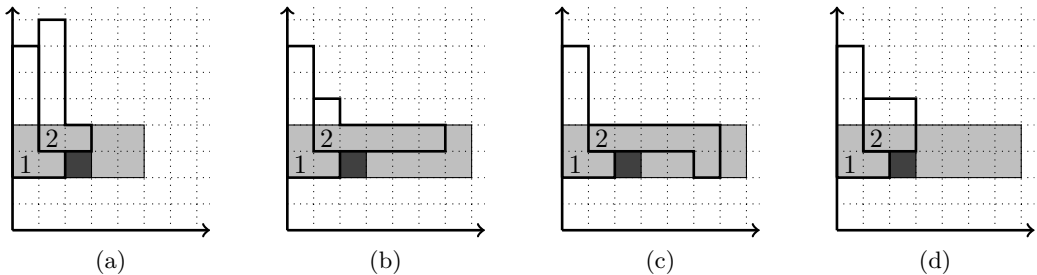

\begin{figure}[h]
\center
\begin{subfigure}{0.2\textwidth}
		\centering
		\begin{tikzpicture}[scale=.35]
        \draw[fill=lightgray] (0,0) rectangle (7,2);
		\draw[line width = 1, ->] (0, -2) -- (7.5, -2);
		\node at (.5, .5) {$1$};
        \node at (1.5, 1.5) {$2$};
        \draw[fill=darkgray] (0,2)--(0,3)--(1,3)--(1,2);
		\draw[line width=1] (0, 0) -- (0, 2) -- (1, 2) -- (1, 1) -- (5, 1) -- (5, 0) -- (0, 0);
        \draw[line width=1] (1,1)--(1,6)--(2,6)--(2,2)--(3,2)--(3,1)--(2,1);
		\draw[dotted,step=1] (0, -2) grid (7.5, 6.5);
		\draw[line width = 1, ->] (0, -2) -- (0, 6.5);
		\end{tikzpicture}
		\caption{}
	\end{subfigure}
~~
\begin{subfigure}{0.2\textwidth}
		\centering
		\begin{tikzpicture}[scale=.35]
        \draw[fill=lightgray] (0,0) rectangle (7,2);
		\draw[line width = 1, ->] (0, -2) -- (7.5, -2);
		\node at (.5, .5) {$1$};
        \node at (1.5, 1.5) {$2$};
        \draw[fill=darkgray] (0,2)--(0,3)--(1,3)--(1,2);
		\draw[line width=1] (0, 0) -- (0,2) -- (1, 2) -- (1, 1) -- (5, 1) -- (5, 0) -- (0, 0);
        \draw[line width=1] (1,3)--(2,3)--(2,2)--(6,2)--(6,1)--(2,1);
		\draw[dotted,step=1] (0, -2) grid (7.5, 6.5);
		\draw[line width = 1, ->] (0, -2) -- (0, 6.5);
		\end{tikzpicture}
		\caption{}
	\end{subfigure}
~~
\begin{subfigure}{0.2\textwidth}
		\centering
		\begin{tikzpicture}[scale=.35]
        \draw[fill=lightgray] (0,0) rectangle (7,2);
		\draw[line width = 1, ->] (0, -2) -- (7.5, -2);
		\node at (.5, .5) {$1$};
        \node at (1.5, 1.5) {$2$};
        \draw[fill=darkgray] (0,2)--(0,3)--(1,3)--(1,2);
		\draw[line width=1] (0, 0) -- (0, 2) -- (1, 2) -- (1, 1) -- (5, 1) -- (5, 0) -- (0, 0);
        \draw[line width=1] (2,1)--(2,6)--(0,6)--(0,5)--(1,5)--(1,1);
		\draw[dotted,step=1] (0, -2) grid (7.5, 6.5);
		\draw[line width = 1, ->] (0, -2) -- (0, 6.5);
		\end{tikzpicture}
		\caption{}
	\end{subfigure}
~
\begin{subfigure}{0.2\textwidth}
		\centering
		\begin{tikzpicture}[scale=.35]
        \draw[fill=lightgray] (0,0) rectangle (7,2);
		\draw[line width = 1, ->] (0, -2) -- (7.5, -2);
		\node at (.5, .5) {$1$};
        \node at (1.5, 1.5) {$2$};
        \draw[fill=darkgray] (0,2)--(0,3)--(1,3)--(1,2);
		\draw[line width=1] (0, 0) -- (0, 2) -- (1, 2) -- (1, 1) -- (5, 1) -- (5, 0) -- (0, 0);
        \draw[line width=1] (1,1)--(3,1)--(3,3)--(1,3)--(1,1);
		\draw[dotted,step=1] (0, -2) grid (7.5, 6.5);
		\draw[line width = 1, ->] (0, -2) -- (0, 6.5);
		\end{tikzpicture}
		\caption{}
	\end{subfigure}
\caption{A $T_3$ tile covers the lower leftmost cell in the gap-base case.}
\label{fig-bis-36}
\end{figure}
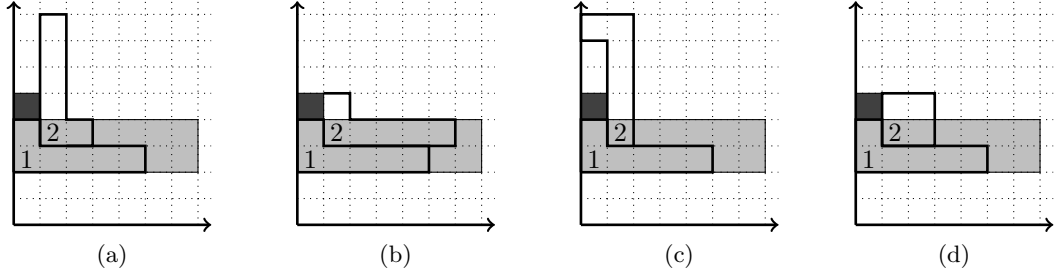

\begin{figure}[h]
\center
\begin{subfigure}{0.3\textwidth}
		\centering
		\begin{tikzpicture}[scale=.35]
        \draw[fill=lightgray] (0,0) rectangle (7,2);
        \draw[fill=lightgray] (-2,2) rectangle (1,4);
		\draw[line width = 1, ->] (-2, -2) -- (7.5, -2);
		\node at (.5, .5) {$1$};
        \node at (1.5, 1.5) {$2$};
		\draw[line width = 1] (-2, 2) -- (0, 2);
        \node at (-.5,3) {\small{GAP}};
		\draw[line width=1] (0, 0) -- (0, 2) -- (1, 2) -- (1, 1) -- (5, 1) -- (5, 0) -- (0, 0);
        \draw[line width=1] (1,1)--(1,6)--(2,6)--(2,2)--(3,2)--(3,1)--(2,1);
		\draw[dotted,step=1] (-2.5, -2.5) grid (7.5, 6.5);
		\draw[line width = 1, ->] (-2, -2) -- (-2, 6.5);
		\end{tikzpicture}
		\caption{}
	\end{subfigure}
~~
\begin{subfigure}{0.3\textwidth}
		\centering
		\begin{tikzpicture}[scale=.35]
        \draw[fill=lightgray] (0,0) rectangle (7,2);
        \draw[fill=lightgray] (-2,2) rectangle (1,4);
		\draw[line width = 1, ->] (-2, -2) -- (7.5, -2);
		\node at (.5, .5) {$1$};
        \node at (1.5, 1.5) {$2$};
		\draw[line width = 1] (-2, 2) -- (0, 2);
        \node at (-.5,3) {\small{GAP}};
		\draw[line width=1] (0, 0) -- (0, 2) -- (1, 2) -- (1, 1) -- (5, 1) -- (5, 0) -- (0, 0);
        \draw[line width=1] (1,1)--(1,5)--(0,5)--(0,6)--(2,6)--(2,1);
		\draw[dotted,step=1] (-2.5, -2.5) grid (7.5, 6.5);
		\draw[line width = 1, ->] (-2, -2) -- (-2, 6.5);
		\end{tikzpicture}
		\caption{}
	\end{subfigure}
~~
\begin{subfigure}{0.3\textwidth}
		\centering
		\begin{tikzpicture}[scale=.35]
        \draw[fill=lightgray] (0,0) rectangle (7,2);
        \draw[fill=lightgray] (-2,2) rectangle (1,4);
		\draw[line width = 1, ->] (-2, -2) -- (7.5, -2);
		\node at (.5, .5) {$1$};
        \node at (1.5, 1.5) {$2$};
		\draw[line width = 1] (-2, 2) -- (0, 2);
        \node at (-.5,3) {\small{GAP}};
		\draw[line width=1] (0, 0) -- (0, 2) -- (1, 2) -- (1, 1) -- (5, 1) -- (5, 0) -- (0, 0);
        \draw[line width=1] (1,1)--(1,3)--(3,3)--(3,1);
		\draw[dotted,step=1] (-2.5, -2.5) grid (7.5, 6.5);
		\draw[line width = 1, ->] (-2, -2) -- (-2, 6.5);
        \draw[line width = 1] (1,3)--(1,4);
		\end{tikzpicture}
		\caption{}
	\end{subfigure}
~~
\begin{subfigure}{0.3\textwidth}
		\centering
		\begin{tikzpicture}[scale=.35]
        \draw[fill=lightgray] (0,0) rectangle (7,2);
		\draw[line width = 1, ->] (-2, -2) -- (7.5, -2);
        \draw[fill=lightgray] (-2,2) rectangle (1,4);
        \node at (-.5,3) {\small{GAP}};
		\node at (.5, .5) {$1$};
        \node at (1.5, 1.5) {$2$};
		\draw[line width=1] (0, 0) -- (0,2) -- (1, 2) -- (1, 1) -- (5, 1) -- (5, 0) -- (0, 0);
        \draw[line width=1] (1,2)--(1,3)--(2,3)--(2,2)--(6,2)--(6,1)--(2,1);
		\draw[dotted,step=1] (-2.5, -2.5) grid (7.5, 6.5);
        \draw[line width = 1] (1,3)--(1,4);
        \draw[line width = 1] (-2,2)--(0,2);
		\draw[line width = 1, ->] (-2, -2) -- (-2, 6.5);
		\end{tikzpicture}
		\caption{}
	\end{subfigure}
~~
\begin{subfigure}{0.3\textwidth}
		\centering
		\begin{tikzpicture}[scale=.35]
        \draw[fill=lightgray] (0,0) rectangle (7,2);
		\draw[line width = 1, ->] (-6, -2) -- (7.5, -2);
       \draw[line width=1] (-6,2)--(0,2);
         \draw[line width=1] (-3,2)--(-3,3)--(2,3)--(2,1);
       \draw[fill=lightgray] (-6,2) rectangle (-3,4);
        \node at (-4.5,3) {\small{GAP}};
		\node at (.5, .5) {$1$};
        \node at (1.5, 1.5) {$2$};
		\draw[line width=1] (0, 0) -- (0, 2) -- (1, 2) -- (1, 1) -- (5, 1) -- (5, 0) -- (0, 0);
		\draw[dotted,step=1] (-6.5, -2.5) grid (7.5, 6.5);
		\draw[line width = 1, ->] (-6, -2) -- (-6, 6.5);
		\end{tikzpicture}
		\caption{}
	\end{subfigure}
\caption{A $T_3$ tile covers the lower leftmost cell in the gap-general case.}
\label{fig-bis-37}
\end{figure}
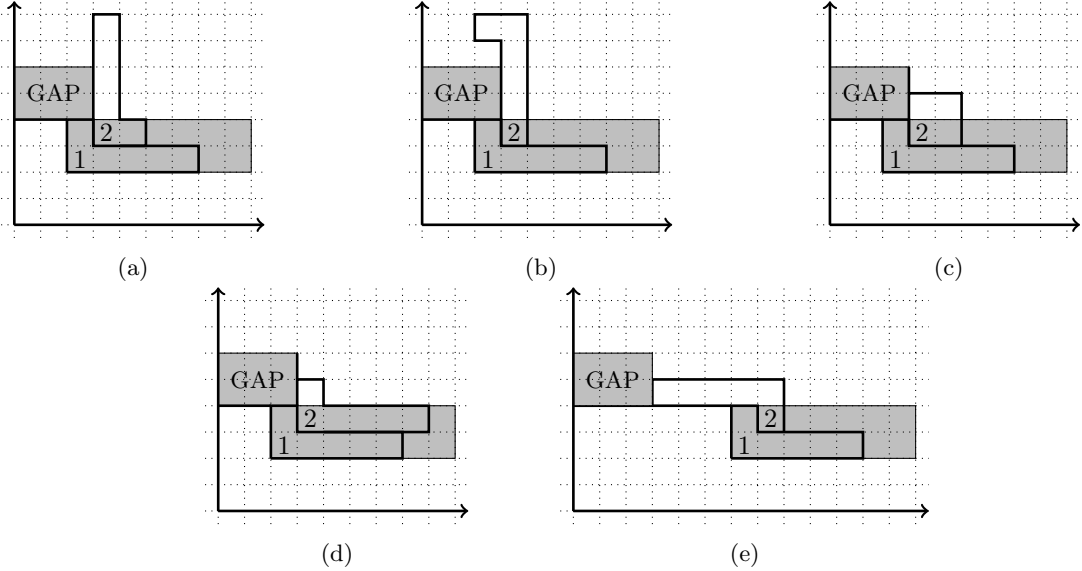

\begin{pf} Let $d$ the distance between the right side of the gap and the $y$-axis. We proceed by induction on $d$. For the induction step, we show that a $T_2$ tile, or a new odd gap of length $L\ge 3$ with the leftmost even $2\times 2$ square not following the rectangular pattern, appears that is above or on the left side of the gap. In the diagrams the gaps are colored in light gray and the cells that cannot be tiled are colored in dark gray.

We consider first the case when the left side of the gap is based on the $y$-axis. This includes the base case $d=3$ for the induction. Look at the tiling of cell $1$, the lowest leftmost in the gap. We cannot use a $T_5$ tile because of the hypothesis. We cannot use a $T_2, T_4$ tile because the gap is too close to the $y$-axis. If we use a $T_1$ tile, then the leftmost even $2\times 2$ square in the gap follows the rectangular pattern. Otherwise, there exists a cell in the lower row of the gap that cannot be tiled. See Figure~\ref{fig-bis-35} for the diagrams of the cases that appear when we try to cover cell $2$, diagonally adjacent to $1$. Note that cell 2 cannot be covered by a $T_2$ tile due to our hypothesis. We only need a right edge of height $1$ for the odd gap. If cell $1$ is tiled by a $T_3$ tile, then cell $2$ cannot be tiled without forcing the leftmost even $2\times 2$ square in the gap to follow the rectangular pattern or leading to a contradiction. See Figure~\ref{fig-bis-36} for the diagrams.
Note that cell 2 cannot be covered by a $T_4$ tile due to our hypothesis.

Consider now the general case. We look at the tiling of cell $1$, the lower leftmost in the gap. We cannot use a $T_5$ tile due to the hypothesis. If we use a $T_2$ tile, the $T_2$ tile is in an irregular position. If we use a $T_1$ tile, then the even $2\times 2$ square containing cell 1 has to follow the rectangular pattern. The reasoning is similar to that done in the case when the left side of the gap is supported by the $y$-axis and the same dark gray cells as in Figure~\ref{fig-bis-35} are impossible to tile. As cell 1 cannot be covered by a $T_4$ tile, it follows that cell 1 is covered by a $T_3$ tile. Consider cell 2 diagonally adjacent to cell 1. See Figure~\ref{fig-bis-37} for the diagrams of the cases that appear. In all cases a new odd gap is created which is closer to the $y$-axis then the original one. If the new odd gap has length 1, then it can be tiled only by a $T_2$ tile in an irregular position. Otherwise, look at the region inside the horizontal strip of width 2 containing the new odd gap that is bounded by the $y$-axis and the left vertical side of the gap. If all even $2\times 2$ squares inside that region follow the rectangular pattern, the new odd gap is forced to have length 1, and a $T_2$ tile in an irregular position appears. If there exists an even $2\times 2$ square that does not follow the rectangular pattern, choose the left side of the new odd gap to be the left side of the leftmost such square. This guarantees that the new odd gap has length $L\ge 3$ and has the leftmost even $2\times 2$ square not following the rectangular pattern.
\end{pf}

\begin{lem}\label{lemma-1234} Assume that the leftmost $2\times 2$-square in an even gap of length $L\ge 2$ does not follow the rectangular pattern. Then there exists a $T_2$ tile in an irregular position that is above the bottom of the gap and to the left of the right side of the gap.
\end{lem}

\begin{pf} If $L=2$, due to the fact that the height of the right side of the gap is at least 2, the tiling of the gap is forced to follow the rectangular pattern. So we may assume $L\ge 4$. If the left side of the gap is on the $y$-axis, a similar analysis to that in the proof of Lemma~\ref{lemma-123} leads to a contradiction. In particular the same dark gray cells marked in~Figure~\ref{fig-bis-35},~Figure~\ref{fig-bis-36} remain impossible to tile.
Consider the general case. We look at the lower leftmost cell in the gap, say 1. Reasoning as in the proof of Lemma \ref{lemma-123}, cell 1 has to be covered by a $T_3$ tile. We look at the tiling of cell 2, diagonally adjacent to cell 1. Following the diagrams in Figure~\ref{fig-bis-37}, this leads to an odd gap that is on the left and above the even gap. Now the existence of the $T_2$ tile follows from Lemma \ref{lemma-123}.
\end{pf}

\begin{lem}\label{lemma9} Assume that a $T_2$ tile is placed in $Q_1$ in an irregular position. Then either all even $2\times 2$ squares to the left and below the $T_2$ tile follow the rectangular pattern, or there exists a $T_2$ tile that is closer to the $y$-axis and is in an irregular position.
\end{lem}

\begin{pf} Figure~\ref{fig:aux-lemma-098} illustrates the statement of the lemma: the dark gray even $2\times 2$ squares follow the rectangular pattern; either the light gray even $2\times 2$ squares follow the rectangular pattern, or there exists a $T_2$ tile in an irregular position that is closer to the $y$-axis.
We proceed by contradiction and assume that not all light gray even $2\times 2$ squares follow the rectangular pattern. We identify the bottom row of width 2 in the light gray region in which appears such a square and apply Lemma \ref{lemma-1234} to an open gap in that row that has a right edge of height 2.
\end{pf}

\begin{lem}\label{lemma-t2} A tiling of $Q_1$ by $\mathcal{T}_n^+$ cannot contain a $T_2$ or $T_4$ tile in an irregular position.
\end{lem}

\begin{pf} Assume that a $T_2$ tile is in an irregular position and at minimal distance from the $y$-axis. By Lemma \ref{lemma9}, all even $2\times 2$ squares to the left and below the $T_2$ tile follow the rectangular pattern. The portion of the top row of the $T_2$ tile that sits between the $T_2$ tile and the $y$-axis has an odd number of cells available. This creates an odd gap to which one applies Lemma \ref{lemma-123} to deduce the existence of a new $T_2$ tile in an irregular position that is based above and to the left of the initial $T_2$ tile. This gives a contradiction.

The statement about the $T_4$ tile follows due to the symmetry of $\mathcal{T}_n$ about the line $y=x$.
\end{pf}

\begin{figure}[h!]
\center
\begin{subfigure}{0.3\textwidth}
		\begin{tikzpicture}[scale=.35]
\draw[fill=lightgray] (-8, 8) rectangle (2, 14);
\draw[fill=gray=6] (-2, 4) rectangle (0, 6);
\draw[fill=gray=6] (0, 4) rectangle (2, 6);
\draw[fill=gray=6] (0, 6) rectangle (2, 8);
\draw[fill=gray=6] (-8, 4) rectangle (2, 8);
\draw[fill=gray=6] (-2, 6) rectangle (0, 8);
\draw[line width=1] (0,8)--(0, 6) -- (2, 6) -- (2, 8) -- (0, 8);
\draw[line width=1] (-8,14)--(2,14);
\draw[line width=1] (-6,12)--(2,12);
\draw[line width=1] (-4,10)--(2,10);
\draw[line width=1] (-6,12)--(-6,14);
\draw[line width=1] (-4,10)--(-4,14);
\draw[line width=1] (-2,8)--(-2,14);
\draw[line width=1] (0,8)--(0,14);
\draw[line width=1] (-2,8)--(-2, 6) -- (0, 6) -- (0, 8) -- (-2, 8);
\draw[line width=1] (-4,8)--(-4, 6) -- (-2, 6) -- (-2, 8) -- (-4, 8);
\draw[line width=1] (-6,10)--(-4, 10) -- (-4, 8) -- (-6, 8);
\draw[line width=1] (-2,6)--(-2, 4) -- (0, 4) -- (0, 6) -- (-2, 6);
\draw[line width=1] (2,4)--(2, 6);
\draw[line width=1] (-6,8)--(-6, 10)--(-8,10)--(-8,12)--(-8,14);
\draw[line width=1] (-8,10)--(-8,4)--(2,4);
\draw[line width=1] (-6,4)--(-6,12);
\draw[line width=1] (-4,4)--(-4,6);
\draw[line width=1] (-8,6)--(-4,6);
\draw[line width=1] (-8,8)--(-4,8);
\draw[line width=1] (-8,12)--(-4,12);
\draw[line width=1] (2,8)--(3,8)--(3,15)--(1,15)--(1,14)--(2,14)--(2,8);
\node at (2.5, 12) {$T_2$};
\draw[dotted] (-8.5, 3.5) grid (5.5, 15.5);
\end{tikzpicture}
		\caption{Lemma \ref{lemma9}.}
		\label{fig:aux-lemma-098}
\end{subfigure}
~
\begin{subfigure}{0.3\textwidth}
\begin{tikzpicture}[scale=0.4]
\draw[dotted] (0, 0) grid (9.5, 8.5);
\draw[line width = 1, ->] (0, 0) -- (9.5, 0);
\draw[line width = 1, ->] (0, 0) -- (0, 8.5);
\node at (-.25, 9.125) {$y$};
\node at (10.25, -.25) {$x$};
\node[below left] at (0, 0) {$(0, 0)$};
\draw[thick] (0, 0) rectangle ++(2, 2);
\draw[thick] (2, 0) rectangle ++(2, 2);
\draw[thick] (0, 2) rectangle ++(2, 2);
\draw[thick] (4, 0) rectangle ++(2, 2);
\draw[thick] (2, 2) rectangle ++(2, 2);
\draw[thick] (0, 4) rectangle ++(2, 2);
\draw[line width = 2.0] (6,0) -- (6,2) -- (4,2) -- (4,4) -- (2,4) -- (2,6) -- (0,6);
\node[below] at (5,0) {};
\draw[dashed] (6,0) rectangle ++(2,2);
\draw[dashed] (4,2) rectangle ++(2,2);
\draw[dashed] (2,4) rectangle ++(2,2);
\draw[dashed] (0,6) rectangle ++(2,2);
\node at (7,1) {$X_1$};
\node at (5,3) {$X_2$};
\node at (3,5) {$X_3$};
\node at (1,7) {$X_4$};
\end{tikzpicture}
\caption{The induction staircase line.}
\label{fig:InductionStep}
\end{subfigure}
\caption{}
\end{figure}

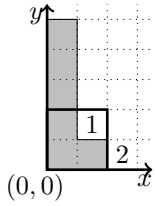
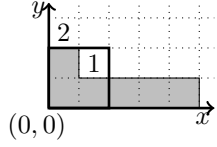
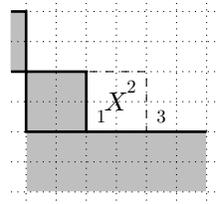
\begin{figure}[h!]
\center
	\begin{subfigure}{0.3\textwidth}
		\centering
		\begin{tikzpicture}[scale=.4]
		\draw[color=white] (4, 1) rectangle (5.5, 5.5);
		\draw[line width = 1., ->] (0, 0) -- (3.5, 0);
		\node at (-.35, 5.125) {$y$};
		\node at (3.25, -.35) {$x$};
		\node at (-.4, -.6) {$(0, 0)$};
		\draw[fill=lightgray] (0, 0) -- (2, 0) -- (2, 1) -- (1, 1) -- (1, 5) -- (0, 5) -- (0, 0);
		\draw[dotted,step=1] (0, 0) grid (3.5, 5.5);
		\draw[line width=1] (0, 0) rectangle ++(2, 2);
		\node at (1.5, 1.5) {$1$};
		\node at (2.5, .5) {$2$};
		\draw[line width = 1, ->] (0, 0) -- (0, 5.5);
		\end{tikzpicture}
		\caption{Base case: $T_1$ covers the corner square.}
		\label{fig:cornercase1-0}
	\end{subfigure}
	~
	\begin{subfigure}{0.3\textwidth}
	\centering
		\begin{tikzpicture}[scale=.4]
		\draw[color=white] (1, 4) rectangle (5.5, 5.5);
		\draw[line width = 1, ->] (0, 0) -- (0,3.5);
		\node at (-.35, 3.25) {$y$};
		\node at (5.125, -.35) {$x$};
		\node at (-.4, -.6) {$(0, 0)$};
		\draw[fill=lightgray] (0, 0) -- (5, 0) -- (5, 1) -- (1, 1) -- (1, 2) -- (0, 2) -- (0, 0);
		\draw[dotted,step=1] (0, 0) grid (5.5, 3.5);
		\draw[line width=1] (0, 0) rectangle ++(2, 2);
		\node at (1.5, 1.5) {$1$};
		\node at (.5, 2.5) {$2$};
		\draw[line width = 1, ->] (0, 0) -- (5.5, 0);
		\end{tikzpicture}
		\caption{Base case: $T_3$ covers the corner square}
		\label{fig:cornercase2-0}
	\end{subfigure}
~
\begin{subfigure}{0.3\textwidth}
\begin{tikzpicture}[scale=.4]
\draw[color=white,fill=lightgray] (0, 0) rectangle ++ (6, 2);
\draw[line width=1] (2,2) -- (6, 2);
\draw[dashed] (2, 2) rectangle ++(2, 2);
\node at (2.5,2.5){\scalebox{0.7}{1}};
\node at (3.5,3.5){\scalebox{0.7}{2}};
\node at (4.5,2.5){\scalebox{0.7}{3}};
\node at (3, 3) {$X$};
\draw[fill=lightgray,color=lightgray] (-.5, 4) rectangle (0, 6);
\draw[line width=1,fill=lightgray] (0, 2) rectangle ++(2, 2);
\draw[line width=1] (-.5, 4) -- (0, 4) -- (0, 6) -- (-.5, 6);
\draw[dotted] (-.5, -.5) grid (6.5, 6.5);
\end{tikzpicture}
\caption{Inductive step}
\label{fig:Illustration}
\end{subfigure}
	\caption{The steps of induction}
	\label{fig:CornerCase}
\end{figure}

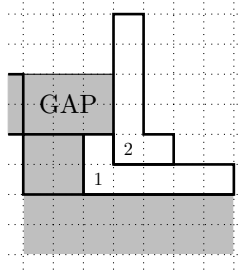
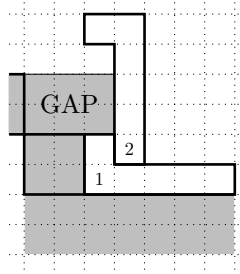
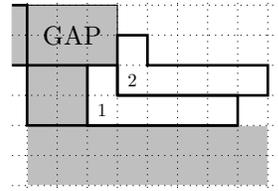
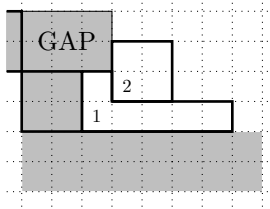
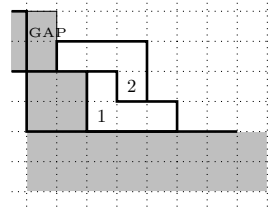
\begin{figure}[h!]
\center
\begin{subfigure}{.2\textwidth}
	\begin{tikzpicture}[scale=.4]
\draw[color=white,fill=lightgray] (0, 0) rectangle ++ (7, 2);
\draw[line width=1] (2,2) -- (7, 2);
\node at (2.5,2.5){\scalebox{0.7}{1}};
\node at (3.5,3.5){\scalebox{0.7}{2}};
\draw[fill=lightgray] (3, 6) rectangle (0, 4);
\node at (1.5, 5) {GAP};
\draw[fill=lightgray,color=lightgray] (-.5, 4) rectangle (0, 6);
\draw[line width=1] (-.5, 4) -- (0, 4) -- (0, 6) -- (-.5, 6);
\draw[line width=1] (2,4)--(3,4)--(3,3)--(7,3)--(7,2)--(2,2)--(2,4);
\draw[line width=1] (3,3)--(5,3)--(5,4)--(4,4)--(4,8)--(3,8)--(3,3);
\draw[line width=1,fill=lightgray] (0, 2) rectangle ++(2, 2);
\draw[dotted] (-.5, -.5) grid (7.5, 8.5);
\end{tikzpicture}
\caption{}
\label{fig:Illustration1}
\end{subfigure}
\qquad\qquad
\begin{subfigure}{.2\textwidth}
\begin{tikzpicture}[scale=.4]
\draw[color=white,fill=lightgray] (0, 0) rectangle ++ (7, 2);
\draw[line width=1] (2,2) -- (7, 2);
\node at (2.5,2.5){\scalebox{0.7}{1}};
\node at (3.5,3.5){\scalebox{0.7}{2}};
\draw[fill=lightgray] (3, 6) rectangle (0, 4);
\node at (1.5, 5) {GAP};
\draw[fill=lightgray,color=lightgray] (-.5, 4) rectangle (0, 6);
\draw[line width=1] (-.5, 4) -- (0, 4) -- (0, 6) -- (-.5, 6);
\draw[line width=1] (2,4)--(3,4)--(3,3)--(7,3)--(7,2)--(2,2)--(2,4);
\draw[line width=1] (3,3)--(4,3)--(4,8)--(2,8)--(2,7)--(3,7)--(3,3);
\draw[line width=1,fill=lightgray] (0, 2) rectangle ++(2, 2);
\draw[dotted] (-.5, -.5) grid (7.5, 8.5);
\end{tikzpicture}
\caption{}
\label{fig:Illustration2}
\end{subfigure}
\qquad\qquad
\begin{subfigure}{.2\textwidth}
\begin{tikzpicture}[scale=.4]
\draw[color=white,fill=lightgray] (0, 0) rectangle ++ (8, 2);
\draw[line width=1] (2,2) -- (7, 2);
\node at (2.5,2.5){\scalebox{0.7}{1}};
\node at (3.5,3.5){\scalebox{0.7}{2}};
\draw[fill=lightgray] (3, 6) rectangle (0, 4);
\node at (1.5, 5) {GAP};
\draw[fill=lightgray,color=lightgray] (-.5, 4) rectangle (0, 6);
\draw[line width=1] (-.5, 4) -- (0, 4) -- (0, 6) -- (-.5, 6);
\draw[line width=1] (2,4)--(3,4)--(3,3)--(7,3)--(7,2)--(2,2)--(2,4);
\draw[line width=1] (3,3)--(8,3)--(8,4)--(4,4)--(4,5)--(3,5)--(3,3);
\draw[line width=1,fill=lightgray] (0, 2) rectangle ++(2, 2);
\draw[dotted] (-.5, -.5) grid (8.5, 6.5);
\end{tikzpicture}
\caption{}
\label{fig:Illustration3}
\end{subfigure}
\qquad\qquad
\begin{subfigure}{.2\textwidth}
\begin{tikzpicture}[scale=.4]
\draw[color=white,fill=lightgray] (0, 0) rectangle ++ (8, 2);
\draw[line width=1] (2,2) -- (7, 2);
\node at (2.5,2.5){\scalebox{0.7}{1}};
\node at (3.5,3.5){\scalebox{0.7}{2}};
\draw[fill=lightgray] (3, 6) rectangle (0, 4);
\node at (1.5, 5) {GAP};
\draw[fill=lightgray,color=lightgray] (-.5, 4) rectangle (0, 6);
\draw[line width=1] (-.5, 4) -- (0, 4) -- (0, 6) -- (-.5, 6);
\draw[line width=1] (2,4)--(3,4)--(3,3)--(7,3)--(7,2)--(2,2)--(2,4);
\draw[line width=1] (3,3)--(5,3)--(5,5)--(3,5)--(3,3);
\draw[line width=1,fill=lightgray] (0, 2) rectangle ++(2, 2);
\draw[dotted] (-.5, -.5) grid (8.5, 6.5);
\end{tikzpicture}
\caption{}
\label{fig:Illustration3-22}
\end{subfigure}
\qquad\qquad
\begin{subfigure}{.2\textwidth}
\begin{tikzpicture}[scale=.4]
\draw[color=white,fill=lightgray] (0, 0) rectangle ++ (8, 2);
\draw[line width=1] (2,2) -- (7, 2);
\node at (2.5,2.5){\scalebox{0.7}{1}};
\node at (3.5,3.5){\scalebox{0.7}{2}};
\draw[fill=lightgray,color=lightgray] (-.5, 4) rectangle (0, 6);
\draw[fill=lightgray] (1, 6) rectangle (0, 4);
\node at (.7, 5.3) {\tiny{GAP}};
\draw[line width=1] (-.5, 4) -- (0, 4) -- (0, 6) -- (-.5, 6);
\draw[line width=1] (2,4)--(3,4)--(3,3)--(5,3)--(5,2)--(2,2)--(2,4);
\draw[line width=1] (4,3)--(4,5)--(1,5)--(1,4);
\draw[line width=1,fill=lightgray] (0, 2) rectangle ++(2, 2);
\draw[dotted] (-.5, -.5) grid (8.5, 6.5);
\end{tikzpicture}
\caption{}
\label{fig:Illustration3-33}
\end{subfigure}
\caption{Inductive step, Case $3$.}
\label{final-ilustration-34}
\end{figure}

\begin{figure}[h!]
\centering
\begin{subfigure}{.3\textwidth}
\begin{tikzpicture}[scale=.4]

\node at (2.5,3) {even$\times$odd};

\node at (2,8) {I};
\node at (5,8) {II};
\node at (7,1) {III};
\node at (7.5,5) {IV};

\draw [line width = 1,->] (0,0)--(0,9);
\draw [line width = 1,->] (0,0)--(9,0);

\draw [line width = 1] (5,2)--(6,2)--(6,7)--(4,7)--(4,6);

\draw [line width = 1] (5,0)--(5,6)--(0,6);

\draw [line width = 1] (4,6)--(4,9);
\draw [line width = 1] (6,6)--(6,9);
\draw [line width = 1] (6,2)--(9,2);
\end{tikzpicture}
\caption{}
\label{final-ilustration-35}
\end{subfigure}
~~~
\begin{subfigure}{.3\textwidth}
\begin{tikzpicture}[scale=.4]

\node at (2,2) {I};
\node at (8,4) {II};

\draw [line width = 1] (0,0)--(10,0)--(10,6)--(0,6)--(0,0);
\draw [line width = 1] (4,0)--(4,4)--(0,4);

\draw [line width = 1] (10,2)--(6,2)--(6,6);

\draw [line width = 1] (1,6)--(1,5)--(5,5)--(5,1)--(9,1)--(9,0);

\draw [line width = 1] (3,4)--(3,5);
\draw [line width = 1] (4,2)--(5,2);
\draw [line width = 1] (5,4)--(6,4);
\draw [line width = 1] (7,1)--(7,2);
\end{tikzpicture}
\caption{}
\label{final-ilustration-39}
\end{subfigure}
\caption{Tilings of $Q_1$ that do not follow the rectangular pattern.}
\end{figure}
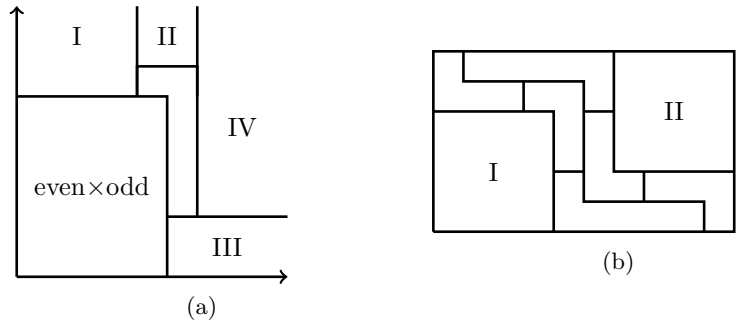

{\it Proof of Theorem \ref{main1}.} We show that every even $2 \times 2$ square follows the rectangular pattern. We do this by induction on a diagonal staircase at $2k$, shown in Figure \ref{fig:InductionStep}. We assume that every even $2 \times 2$ square southwest of this line satisfies the hypothesis and we prove that every even $2\times 2$ square $X_i$ in Figure \ref{fig:InductionStep} also satisfies it. We first investigate the tiling of the corner cell of $Q_1$. If tiled by a $T_5$ tile, we are done. The other possible cases are shown in Figures \ref{fig:cornercase1-0}, \ref{fig:cornercase2-0}. Assume $T_1$ covers the corner square. If cell $1$ is covered by $T_1, T_3,T_4$ or $T_5$, then cell $2$ cannot be covered by any tile in $\mathcal{T}_n^+$. Thus cell $1$ has to be covered by $T_2$, and we  complete an $n\times 2$ rectangle. The other case in Figure \ref{fig:cornercase2-0} is solved similarly, via a symmetry about $x=y$.

For the induction step we prove that the even $2\times 2$ squares $X_i$ in Figure \ref{fig:InductionStep} follow the rectangular pattern. Choose the rightmost square $X_i,$ say $X$, which does not follow the rectangular pattern (see Figure \ref{fig:Illustration}).  Note that $X$ is bounded below and two units to the right by the $x$ axis or by two even $2\times 2$ squares that follow the rectangular pattern, and it is bounded to the left by the $y$-axis or an even $2\times 2$ square that follows the rectangular patter. By assumption, cell $1$ cannot be tiled by a $T_5$ tile. Also, cell $1$ cannot be tiled by $T_4$ due to Lemma~\ref{lemma-t2}, as the $T_4$ tile is in irregular position. We discuss the other cases below.

\emph{Case 1.} Assume that a $T_1$ tile covers cell $1$ in Figure \ref{fig:Illustration}. If cell $2$ is covered by $T_1,T_3,T_4,$ or $T_5$, then cell 3 is impossible to cover by $\mathcal{T}_n^+$. If cell $2$ is covered by $T_2$, then the square $X$ is covered by an $n\times 2$ rectangle.

\emph{Case 2.} If a $T_2$ tile covers cell $1$ in Figure \ref{fig:Illustration}, apply Lemma~\ref{lemma-t2}.

\emph{Case 3.} Assume that a $T_3$ tile covers cell $1$ in Figure \ref{fig:Illustration}.

\emph{Subcase 1.} Assume that a $T_1, T_2, T_3,$ or $T_5$ tile covers cell $2$. We work now with Figure~\ref{final-ilustration-34}. Then an odd gap appears on the left side of the tile covering cell 2. From Lemma~\ref{lemma-123} it follows that there exists a $T_2$ tile  in irregular position, which by Lemma~\ref{lemma-t2} is in contradiction to the existence of a tiling for the $Q_1$.

\emph{Subcase 2.} Assume that a $T_4$ tile covers cell $2$. This completes a $2\times n$ rectangle that covers $X$.

Finally, we show that the addition of an extra even $\times$ odd or odd $\times$ odd rectangle to $\mathcal{T}_n^+$ or $\mathcal{T}_n$ allows for a tiling of $Q_1$ that does not respect the rectangular pattern. As the concatenation of two odd $\times$ odd rectangles is an even $\times$ odd rectangle, we can consider only the last type. Also, a concatenation of an odd number of copies of an even $\times$ odd rectangles can be used to construct an even $\times$ odd rectangle of arbitrary large length and height. Assuming the existence of such a rectangle in the tiling, a tiling of $Q_1$ that does not follow the rectangular pattern is shown in Figure~\ref{final-ilustration-35}. The base of the even $\times$ odd rectangle is odd and the height is even. The regions I, II, III are half-infinite strips of even width and region IV is a copy of the $Q_1$. All of them can be tiled by $2\times n$ or $n\times 2$ rectangles.

\bibliographystyle{plain}

\end{document}